\documentclass[12pt,a4paper]{article}
\usepackage{float}
\usepackage{amsmath}
\usepackage{amsfonts}
\usepackage{amssymb}
\usepackage{graphicx}
\renewcommand{\baselinestretch}{1}
\newtheorem{theorem}{Theorem}

\newtheorem{lemma}{Lemma}

\newenvironment{proof}[1][Proof]{\noindent\textbf{#1.} }{\ \rule{0.5em}{0.5em}}

\let \a = \alpha
\let \b = \beta
\let \l = \lambda
\let \g = \gamma

\begin{document}

\title{On the concave and convex solutions of mixed convection boundary layer approximation in a porous medium}
\author{BERNARD BRIGHI$\dag$  and JEAN-DAVID HOERNEL$\ddagger$ }
\date{}
\maketitle

\begin{center}
Universit\'e de Haute-Alsace, Laboratoire de Math\'ematiques et Applications
\vskip 0,1cm
4 rue des fr\`eres Lumi\`ere, 68093 MULHOUSE (France)
\vskip 1cm
\end{center}

\begin{abstract}
In this paper we investigate the similarity solutions of a plane mixed convection boundary layer flow near a semi-vertical plate, with a prescribed power law function of the distance from the leading edge for the temperature, that is embedded in a porous medium. We show the existence and uniqueness of convex and concave solutions for positive values of the power law exponent.
\end{abstract}

\renewcommand{\baselinestretch}{1}
\footnotetext{\hspace{-0.8cm} AMS 2000 Subject Classification: 34B15, 34C11, 76D10.}
\footnotetext{\hspace{-0.8cm} Key words and phrases: Boundary layer, similarity solution, third order nonlinear differential equation,  boundary value problem. }
\footnotetext{\hspace{-0.8cm} $\dag$ bernard.brighi@uha.fr $\ddagger$ j-d.hoernel@wanadoo.fr}

\section{Introduction}
We look at the nonlinear autonomous differential equation established in \cite{aly}
\begin{equation}
f'''+(1+\lambda)ff''+2\lambda(1-f')f'=0 \label{eq}
\end{equation}
on $[0,\infty)$ with the boundary conditions
\begin{equation}
f(0)=\alpha, \quad f'(0)=\beta \quad \text{and}\quad f'(\infty)=1 \label{c1}
\end{equation}
where $\b>0$ and $f'(\infty):=\underset{t\rightarrow\infty}{\lim} f'(t)$.

The problem (\ref{eq})-(\ref{c1}) with $\a=0$ and $\b=1+\varepsilon$, $-1<\varepsilon<\frac{1}{2}$ is considered in \cite{aly} and comes from the study of the mixed convection boundary-layer flow along a semi-infinite vertical plate embedded in a saturated porous medium. We take a rectangular Cartesian co-ordinate system with the origin fixed at the leading edge of the vertical plate, the $x$-axis directed upward along the plate and the $y$-axis normal to it and we write the temperature on the plate as $T(x,0)=T_\infty+Ax^\lambda$ with $A>0$, $\lambda\in \mathbb R$ and where $T_\infty$ denotes the tempetature far from the wall. The fluid velocity at the edge of the boundary layer is writen as $Bx^\lambda$ with $B>0$ and we suppose that there is no normal velocity on the plate. Then, for large P\'eclet numbers, one can compute the stream function and the temperature of the fluid in the porous medium using the dimensionless similarity function $f$ given by (\ref{eq})-(\ref{c1}). For more details on the derivation of the model and numerical results, see \cite{aly} and \cite{nazar}.

The solution of (\ref{eq})-(\ref{c1}) depends on two parameters, the power law exponent $\lambda$ and the mixed convection parameter $\varepsilon=\frac{Ra}{Pe}$ with $Ra$ the Rayleigh number and $Pe$ the P\'eclet number. For $\varepsilon=0$ we have forced convection, $\varepsilon>0$ corresponds to aiding mixed convection and $\varepsilon<0$ to opposing mixed convection.

For $\lambda=0$ (isothermal surface), equation (\ref{eq}) reduces to the Blasius equation (see \cite{bla}) and is studied with conditions (\ref{c1}) in \cite{coppel} for $\alpha\geq 0$ and $\beta>0$, Êin \cite{hart} for $\alpha\in\mathbb R$ and $\beta\in(0,1)$ and in \cite{brighi03} for $\alpha\in\mathbb R$ and $\beta>1$. Blasius equation is also a particular case of the Falkner-Skan equation $f'''+ff''+\lambda(1-f'^2)=0$ (see \cite{hart}).

The equivalent problem of free convection over a vertical flate plate embedded in a fluid saturated porous medium  is studied in \cite{ing}, \cite{brighi02}, \cite{brighi01} and \cite{BrighiSari}.

In \cite{guedda1} one can found interesting new results about the problem (\ref{eq})-(\ref{c1}) for $\lambda<0$ and $\a=0$, in particular that there are infinitely many solutions for $-1<\lambda<0$ and $0<\beta<\frac{3}{2}$, no solutions for $\lambda =-1$ (except the trivial one $f(t)=t$ for $\b=1$) and no nonnegative solutions for $\lambda <-1$ and $\beta>\frac{3}{2}$.

In this note, we are interested in the case $\lambda>0$, and prove existence and uniqueness of the solution of (\ref{eq})-(\ref{c1}) whose second derivative does not vanish.

\section{Existence and uniqueness results}

First, let us suppose that $f$ satisfies equation (\ref{eq}) on some interval $I$ and denote by $F$ any anti-derivative of $f$ on $I$, then we have
\begin{equation}
\left ( f''e^{(1+\lambda)F}\right)'=-2\lambda(1-f')f'e^{(1+\lambda)F}. \label{exp1}
\end{equation}
As for the Falkner-Skan equation, and contrary to what occurs in \cite{brighi01}, the relation (\ref{exp1}) shows that oscillatory solutions may exist (see \cite{hastings}). Nevertheless, for the problem (\ref{eq})-(\ref{c1}) we could expect convex solutions (i.e. such that $f''>0$) if $0<\b<1$, and concave solutions (i.e. such that $f''<0$) if $\b>1$. Before proving this, let us look briefly at the case $\lambda>0$ and $\beta=1$; the function $g(t)=t+\alpha$ is then a solution of (\ref{eq})-(\ref{c1}).
Let $f$ be another solution with $f''(0)=\gamma$ and assume first that $\gamma>0$. Then, since $f'(0)=f'(\infty)=1$, there exists $t_0>0$ such that $f'(t_0)>1$, $f''(t_0)=0$ and $f'''(t_0)\leq 0$. However, from (\ref{eq}) we obtain $f'''(t_0)=-\lambda f'(t_0)(1-f'(t_0))>0$ and thus a contradiction.
If $\gamma<0$, the same approach does not allow us to produce a contradiction, but leads to $f'(t_0)<0$, and thus $g$ is the unique increasing solution of (\ref{eq})-(\ref{c1}) when $\lambda>0$ and $\beta=1$. Clearly, at this stage, we cannot exclude that oscillatory solutions exist.

In order to obtain solutions of (\ref{eq})-(\ref{c1}), we consider the following initial value problem:
\begin{equation}
\left \{ \begin{array}
[c]{lll}
f'''+(1+\lambda)ff''+2\lambda(1-f')f'=0,\\
f(0) = \a,\\
f'(0) = \b,\\
f''(0) = \gamma
\end{array} \right . \label{initial}
\end{equation}
with $\lambda>0$, $\alpha\in\mathbb R$, $\beta>0$ and a suitable $\gamma$. We denote the solution by $f_\gamma$ and by $[0,T_\gamma)$ its right maximal interval of existence. Next, integrating (\ref{eq}) on $[0,t]$ for $0<t<T_\gamma$, we obtain the useful identity
\begin{equation}
f''_\gamma(t)-\gamma+(1+\lambda)(f_\gamma(t)f'_\gamma(t)-\alpha\beta)
+2\lambda(f_\gamma(t)-\alpha)=(3\lambda+1)\int_0^t f'_\gamma(s)^2ds. \label{i1}
\end{equation}
We also need the following lemma.

\begin{lemma}\label{vanish}
If $f$ is a solution of {\rm (\ref{eq})} on $[0,T_\gamma)$ such that there exists a point $t_0$ satisfying $f''(t_0)=0$ and $f'(t_0)=1$, then $f''(t)=0$ for every $t\in [0,T_\gamma)$.
\end{lemma}
\begin{proof}
Let $f$ be a solution of (\ref{eq}) on $[0,T_\gamma)$ such that $f''(t_0)=0$ and Ê$f'(t_0)=1$ for some $t_0\in[0,T_\gamma)$. Since the function $g(t)=t-t_0+f(t_0)$ is a solution of (\ref{eq}) such that $g(t_0)=f(t_0)$, $g'(t_0)=f'(t_0)$ and $g''(t_0)=f''(t_0)$, we obtain $g=f$ and $f''\equiv 0$.
\end{proof}
\begin{theorem}
Let $\lambda>0$. For $\a\in\mathbb R$ and $0<\b<1$ problem {\rm (\ref{eq})-(\ref{c1})} admits a unique convex solution.
\end{theorem}
\begin{proof}[Proof of existence] Let $f_\gamma$ be a solution of the initial value problem (\ref{initial}) with $0<\b<1$ and $\gamma\geq 0$. We notice that $f_\gamma(t)$ exists as long as we have $f''_\g>0$ and $f'_\g<1$. From Lemma \ref{vanish}, $f''_\gamma$ cannot vanish at a point where $f'_\gamma=1$. Therefore, it follows that there are only three possibilities:

\vspace{3mm}
(a) $f''_\gamma$ becomes negative from a point such that $f'_\gamma<1,$

(b) $f'_\gamma$ takes the value $1$ at some point for which $f''_\gamma>0,$

(c) we always have $0<f'_\gamma<1$ and $f''_\gamma>0.$
\vspace{3mm}

\noindent As $f'_0(0)=\b<1$, $f''_0(0)=0$ and $f'''_0(0)=2\lambda \b(\b-1)<0$, we have that $f_0$ is of type (a), and by continuity it must be so for $f_\gamma$ with $\gamma>0$ small enough.

On the other hand, as long as $f''_\gamma(t)>0$ and $f'_\gamma(t)\leq 1$, we have  $f_\gamma(t)\leq t+\a$, and (\ref{i1}) leads to
\begin{align*}
f''_\gamma(t) &\geq -(1+\lambda)f_\gamma(t)f'_\gamma(t)-2\lambda f_\gamma(t)+(1+\lambda)\a\b+2\lambda \a+\gamma \cr
&\geq-(3\lambda+1)t +(1+\lambda)(\a\b-|\a|)+\gamma
\end{align*}
and integrating once again we have
\begin{equation}
f'_\gamma(t) \geq -\frac{3\lambda+1}{2}t^2+((1+\lambda)(\a\b-|\a|)+\gamma)t+\b.\label{p}
\end{equation}
Hence, for $\gamma$ large enough, the polynomial on the right of (\ref{p}) takes values greater than 1. Therefore, for such a $\g$, there exists $t_0$ such that $f'_\gamma(t_0)=1$ and $f''_\gamma(t)>0$ for $t\leq t_0$, and $f_\gamma$ is of type (b).

Defining $A=\left \{ \gamma>0\, ; \, \text{$f_\gamma$ is of type (a)}\right \}$ and
$B=\left \{ \gamma>0\, ; \, \text{$f_\gamma$ is of type (b)}\right \}$ we have that
$A \neq \emptyset$, $B \neq \emptyset$ and $A \cap B = \emptyset$. Both $A$ and $B$ are open sets, so there exists a $\gamma_*>0$ such that the solution $f_{\gamma_*}$ of
(\ref{initial}) is of type (c) and is defined on the whole interval $[0,\infty)$. For this solution we have that $0<f'_{\gamma_*}<1$ and $f''_{\gamma_*}>0$ which implies that $f'_{\gamma_*}(t) \rightarrow l\in(0,1]$ as $t\rightarrow \infty$. Suppose that $l\neq1$, as $f'_{\gamma_*}$ is increasing we have $f'_{\gamma_*}\leq l$ and
\begin{align*}
f'''_{\gamma_*}(t)&\leq -(1+\lambda)\a f''_{\gamma_*}(t)-2\lambda(1-l)f'_{\gamma_*}(t).
\end{align*}
Integrating this inequality leads to
\begin{align*}
f''_{\gamma_*}(t))&\leq  \gamma_*-(1+\lambda)\a(f'_{\gamma_*}(t)-\b)-2\lambda(1-l)(f_{\gamma_*}(t)-\a)
\end{align*}
and, as $f'_{\gamma_*}(t) \rightarrow l<1$ and  $f_{\gamma_*}(t)\rightarrow \infty$ as $t\rightarrow \infty$, we obtain a contradiction with the positivity of $f''_{\gamma_*}$.
\end{proof}

\vspace{3mm}
\begin{proof}[Proof of uniqueness]
Let $f$ be a convex solution of (\ref{eq})-(\ref{c1}). As $f(t)>0$ for $t$ large enough, we have that $f'''(t)<0$ for $t$ large enough. Then $f''(t)\rightarrow 0$ as $t\rightarrow \infty$.

\noindent As $f'$ and $f''$ are positive, we can define a function $v:[\b^2,1)\rightarrow[\a,\infty)$ such that
$$\forall t\geq 0,\quad v(f'(t)^2)=f(t).$$
Setting $y=f'(t)^2$ leads to
\begin{equation}
f(t)=v(y),\quad f''(t)=\frac{1}{2v'(y)}\quad \text{and} \quad f'''(t)=-\frac{v''(y)\sqrt{y}}{2v'(y)^3}\label{y}
\end{equation}
and using (\ref{eq}) we obtain
\begin{equation}
\forall y\in [\b^2,1),\quad v''=(1+\l)\frac{vv'^2}{\sqrt{y}}+4\l\left (1-\sqrt{y} \right)v'^3 \label{eq-v}
\end{equation}
with
$$v(\b^2)=v(f'(0)^2)=\a,\quad v(1)=\underset{t\rightarrow \infty}{\lim}\, f(t)=\infty  \quad \text{and}\quad v'(\b^2)=\frac{1}{2\g}>0.$$
Suppose now that there are two convex solutions $f_{1}$ and $f_{2}$ of
(\ref{eq})-(\ref{c1}) with $f''_i(0)=\gamma_i>0$, $i\in\{1,2\}$ and $\gamma_1>\gamma_2$. They gives $v_{1},\, v_{2}$ solutions of equation (\ref{eq-v}) defined on
$[\b^2,1)$ such that
$$v_{1}(\b^2)=v_{2}(\b^2)=\a, \quad v'_{1}(\b^2)=\frac{1}{2\gamma_1}  \quad \text{and}\quad v'_{2}(\b^2)=\frac{1}{2\gamma_2}.$$
Let $w=v_{1}-v_{2}$, we have $w(\b^2)=0$ and $w'(\b^2)< 0$.
If $w'$ vanishes, there exists an $x$ in $[\b^2,1)$ such that $w'(x)=0$, $w''(x)\geq0$ and $w(x)<0$. But from (\ref{eq-v}) we then obtain
$$w''(x)=(1+\l)\frac{v'_1(x)^2}{\sqrt{x}}w(x)<0$$
and this is a contradiction. Therefore $w'<0$ and $w<0$ on $[\b^2,1)$. Set now $V_i=\frac{1}{v'_i}$ for $i\in\{1,2\}$ and $W=V_1-V_2$. We have $W>0$ and
using (\ref{eq-v}) we obtain
$$W'(y)=-(1+\l)\frac{w(y)}{\sqrt{y}}-4\l(1-\sqrt{y})w'(y)>0.$$
But, using (\ref{y}) we have $V_i(f'(t)^2)=2f''(t)$ and thus $V_i(y)\rightarrow 0$ as $y\rightarrow 1$. Hence $W(y) \rightarrow 0$ as $y\rightarrow 1$, a contradiction.
\end{proof}
\begin{theorem}
Let $\lambda>0$. For $\a \in \mathbb R$ and $\b>1$, problem {\rm (\ref{eq})-(\ref{c1})} admits a unique  concave solution.
\end{theorem}
\begin{proof}[Proof of existence]
Let $f_\gamma$ be a solution of the initial value problem (\ref{initial}) with $\b>1$ and $\gamma\leq 0$. As long as we have $f'_\gamma>1$ and $f''_\gamma<0$, then $f_\gamma$ exists. Because of Lemma \ref{vanish}, there are only three possibilities

\vspace{3mm}
(a) $f''_\gamma$ becomes positive from a point such that $f'_\gamma>1,$

(b) $f'_\gamma$ takes the value $1$ at some point for which $f''_\gamma<0,$

(c) we always have $1<f'_\gamma$ and $f''_\gamma<0.$

\vspace{3mm}
\noindent As $f'_0(0)=\b>1$, $f''_0(0)=0$ and $f'''_0(0)=2\lambda \b(\b-1)>0$ we have that $f'_0(t)>1$ and $f''_0(t)>0$ on some interval $[0,t_0)$. Then, by continuity for small values of $-\gamma$ we have that $f''_\gamma$ becomes positive at some point with $f'_\gamma>1$ and $f_\gamma$ is of type (a).

As long as $f''_\gamma(t)<0$ and $f'_\gamma(t)\geq 1$, we have $f_\gamma(t)\geq \a$ and (\ref{i1}) leads to
\begin{equation*}
f''_\gamma(t) \leq \g +(1+\l)\a\b+(3\l+1)\b^2t+|\a|(1+\l)\b.
\end{equation*}
Integrating once again we have
$$f'_\gamma(t) \leq \frac{3\lambda+1}{2}\b^2t^2+(\g+(1+\lambda)\a\b+|\a|(1+\l)\b)t+\b.$$
Hence, for $-\gamma$ large enough, there exists $t_0$ such that $f'_\gamma(t_0)=1$ and $f''_\gamma(t)<0$ for $t\leq t_0$, and $f_\gamma$ is of type (b).

Defining $A=\left \{ \gamma<0\, ; \, \text{$f_\gamma$ is of type (a)}\right \}$ and
$B=\left \{ \gamma<0\, ; \, \text{$f_\gamma$ is of type (b)}\right \}$ we have that
$A \neq \emptyset$, $B \neq \emptyset$ and $A \cap B = \emptyset$. Both $A$ and $B$ are open sets, so there exists a $\gamma_*<0$ such that the solution $f_{\gamma_*}$ of (\ref{initial}) is of type (c) and is defined on the whole interval $[0,\infty)$. For this solution we have that $f'_{\gamma_*}>1$ and $f''_{\gamma_*}<0$ which implies that $f'_{\gamma_*} \rightarrow l\geq 1$ as $t\rightarrow \infty$. Suppose that $l\neq1$, as $f'_{\gamma_*}$ is decreasing we have $f'_{\gamma_*}\geq l$ and
\begin{align*}
f'''_{\gamma_*}(t)&\geq -(1+\lambda)\a f''_{\gamma_*}(t)-2\lambda(1-f'_{\gamma_*}(t))f'_{\gamma_*}(t).
\end{align*}
Integrating this inequality leads to
$$f''_{\gamma_*}(t)\geq  \gamma_*-(1+\lambda)\a(f'_{\gamma_*}(t)-\b)-2\lambda(1-l)(f_{\gamma_*}(t)-\a)$$
and, as $f'_{\gamma_*}(t) \rightarrow l>1$ and  $f_{\gamma_*}(t)\rightarrow \infty$ as $t\rightarrow \infty$, we obtain a contradiction with the negativity of $f''_{\gamma_*}$.
\end{proof}

\vspace{3mm}
\begin{proof}[Proof of uniqueness]
Let $f_1$ and $f_2$ be two concave solutions of {\rm (\ref{eq})-(\ref{c1})} and let $\gamma_i=f''_i(0)<0$, $i\in \left\{ 1,2 \right \}$ with $\gamma_1>\gamma_2$. Writing $g=f_1-f_2$, we have $g'(0)=0$, $g'(\infty)=0$ and $g''(0)>0.$
Hence there exists $t_0>0$ such that $g'$ admits a positive maximum at $t_0$, and this means
$$g'(t_0)>0, \, g''(t_0)=0 \text{ and } g'''(t_0)\leq 0.$$
From (\ref{eq}), and since $f'_i>1$ and $f''_i<0$, we obtain
$$g'''(t_0)=-(1+\lambda)f''_1(t_0)g(t_0)-2\lambda g'(t_0)(1-f'_1(t_0)-f'_2(t_0))>0$$
and hence a contradiction.
\end{proof}

\end{document}